\documentclass[leqno,12 pt]{article}

\usepackage{amssymb,latexsym}
\usepackage{amsmath}
\usepackage{amsfonts}
\usepackage{amsthm}
\usepackage{amssymb}
\usepackage{color}
\usepackage{verbatim}
\usepackage{comment}
\usepackage{graphicx}
\usepackage[margin=3cm]{geometry}
\usepackage{hyperref}
\newcommand{\footremember}[2]{%
    \footnote{#2}
    \newcounter{#1}
    \setcounter{#1}{\value{footnote}}%
}
\newcommand{\footrecall}[1]{%
    \footnotemark[\value{#1}]%
}

\usepackage{indentfirst}

\numberwithin{equation}{section}

\theoremstyle{plain}
\newtheorem{theorem}{Theorem}
\newtheorem{lemma}[theorem]{Lemma}
\newtheorem{proposition}[theorem]{Proposition}
\newtheorem{corollary}[theorem]{Corollary}

\newtheorem{question}{Question}

\theoremstyle{definition}

\newtheorem{remark}[equation]{Remark}

\def\N{\mathbb{N}}
\def\Z{\mathbb{Z}}
\def\Q{\mathbb{Q}}

\title{Triples which are $D(n)$-sets for several n's}

\begin{document}

\date{}


\author{%
Nikola Ad\v{z}aga\footremember{alley}{Faculty of Civil Engineering, University of Zagreb; nadzaga@grad.hr}%
\and Andrej Dujella \footremember{trailer}{Department of Mathematics, Faculty of Science, University of Zagreb; duje@math.hr}%
\and Dijana Kreso \footrecall{alley} \footnote{Department of Mathematics, University of British Columbia; kreso@math.ubc.ca}%
\and Petra Tadi\'c\footrecall{trailer} \footnote{Department of Mathematics, Statistics and Information Science, Juraj Dobrila University of Pula; ptadic@unipu.hr}%
}


\maketitle
\begin{abstract}
 For a nonzero integer $n$, a set of distinct nonzero integers $\{a_1,a_2,\ldots,a_m\}$ such that  $a_ia_j+n$ is a perfect square for all $1\leq i<j\leq m$, is called a  Diophantine $m$-tuple with the property $D(n)$ or simply $D(n)$-set. $D(1)$-sets are known as  simply Diophantine $m$-tuples. Such sets were first studied by Diophantus of Alexandria, and since then by many authors.  It is natural to ask if there exists a Diophantine $m$-tuple (i.e.\@ $D(1)$-set) which is also a $D(n)$-set for some $n\neq 1$. This question was raised by Kihel and Kihel in 2001. They conjectured that  there are no Diophantine triples which are also $D(n)$-sets for some $n\neq 1$. However, the conjecture does not hold, since, for example, $\{8, 21, 55\}$ is a $D(1)$ and $D(4321)$-triple, while $\{1, 8, 120\}$ is a $D(1)$ and $D(721)$-triple. 
We present several infinite families of Diophantine triples  $\{a, b, c\}$ which are also  $D(n)$-sets for two distinct $n$'s with $n\neq 1$, as well as some Diophantine triples which are also $D(n)$-sets for three distinct $n$'s with $n\neq 1$. We further consider some related questions.
\end{abstract}



\section{Introduction}
A set of distinct nonzero integers $\{a_1,a_2,\ldots,a_m\}$ such that  $a_ia_j+1$ is a perfect square for
all $1\leq i<j\leq m$, is called a \emph{Diophantine $m$-tuple}.
 The problem of constructing such sets was first studied by
Diophantus of Alexandria, who found four rational numbers with the same property. Fermat found a first Diophantine quadruple --  the set $\{1,3,8,120\}$. Any Diophantine pair $\{a,b\}$ can be extended to a Diophantine triple, e.g.\@ by adding $a + b + 2r$ to the set, where $ab+1=r^2$.  Also, any Diophantine triple $\{a,b, c\}$ can be extended to a Diophantine quadruple. Namely, let $ab+1=r^2$, $bc+1=s^2$, $ca+1=t^2$, where $r, s, t$ are positive integers. Then for $d=a+b+c+2abc+2rst$, the set $\{a, b, c, d\}$ is a Diophantine quadruple. Quadruples of this form are called \emph{regular}. In 2004, Dujella~\cite{d05} showed that there are no Diophantine sextuples and that there are at most finitely many Diophantine quintuples.  A couple of decades old conjecture that there are no Diophantine quintuples seems to have been confirmed, according to a recent preprint~\cite{HTZ16+}.  A stronger and still open conjecture is that all Diophantine quadruples are regular.

Various generalizations of the problem have been studied in the literature. For a nonzero integer $n$, a set of distinct nonzero integers $\{a_1,a_2,\ldots,a_m\}$ such that  $a_ia_j+n$ is a perfect square for all $1\leq i<j\leq m$, is called a  Diophantine $m$-tuple with the property $D(n)$ or simply $D(n)$-set. Note that a Diophantine $m$-tuple is a $D(1)$-set. 

The best known upper bounds for the size of  $D(n)$-sets are logarithmic in $n$~\cite{dsize,dsize2}. If $n$ is a prime, it is known that there is an absolute bound for the size of a $D(n)$-set~\cite{dl}. For an arbitrary $n\equiv 2 \pmod 4$ it is easy to see that that there does not exist a $D(n)$-quadruple, while if $n\not \equiv 2 \pmod 4$ and $n\not \in S=\{-4, -3, -1, 3, 5, 12, 20\}$, there exists at least one $D(n)$-quadruple \cite{dacta}. It is conjectured that for $n\in S$ there does not exist a $D(n)$-quadruple.

It is natural to ask if there exists a Diophantine $m$-tuple (i.e. $D(1)$-set) which is also a $D(n)$-set for some $n\neq 1$. This question was raised by Kihel and Kihel in \cite{kk01}. They conjectured that  there are no Diophantine triples which are also $D(n)$-sets for some $n\neq 1$. However, the conjecture does not hold, since, for example, $\{8, 21, 55\}$ is a $D(1)$ and $D(4321)$-triple (as noted in the MathSciNet review
of \cite{kk01}), while $\{1, 8, 120\}$ is a $D(1)$ and $D(721)$-triple. The latter example was found by Zhang and Grossman~\cite{zg15}.

Now, one may ask for how many different $n$'s  with $n\neq 1$  can a $D(1)$-set  also be a $D(n)$-set. One easily sees that there are only finitely many such $n$'s for a fixed $D(1)$-set $\{a, b, c\}$. One way to see that is to note that on the elliptic curve
\begin{equation}\label{curve}
E(\Q)\,: \quad y^2=(x+ab)(x+ac)(x+bc),
\end{equation}
for which we say that it is \emph{induced} by a Diophantine triple $\{a, b, c\}$, there are only finitely many integer points. An elementary proof of this fact follows from the results of \cite{zg15}.

Here we present several infinite families of Diophantine triples (i.e.\ $D(1)$-sets) which are also  $D(n)$-sets for two distinct $n$'s with $n\neq 1$ (Section~\ref{sec2}).
We further show that there are Diophantine triples which are also $D(n)$-sets for three distinct $n$'s with $n\neq 1$ (Section~\ref{sec3}). We do so by reducing the problem to a problem about induced elliptic curves.  Note that for a Diophantine triple $\{a, b, c\}$,  the curve $E(\Q)$ has several obvious rational points, e.g.\@ $P=(0, abc)$ and $S=(1, rst)$, where $r, s, t$ are defined by $bc+1=r^2$, $ca+1=s^2$, $ab+1=t^2$. Furthermore, one can show that $S\in 2E(\Q)$, and that for $T\in 2E(\Q)\cap \Z^2$, $\{a, b, c\}$ is a  $D(x(T))$-set. See Section~\ref{sec2} for details. In our search for Diophantine triples which are also $D(n)$-sets for several distinct $n$'s with $n\neq 1$, we thus focus on Diophantine triples $\{a, b,c\}$ for which  $2kP+\ell S\in \Z^2$ with $k,\ell \in \Z$.

Motivated by the fact that the size of a set $N$ for which there exists a triple $\{a, b, c\}$ of nonzero integers (not necessarily a Diophantine triple)
which is a $D(n)$-set for all $n\in N$ can be arbitrarily large, and our unsuccessful search for a Diophantine triple which is also a $D(n)$-set for four distinct $n$'s with $n\neq 1$, in the last section we consider the following question and present some computational results: For a given positive integer $k$, what can be said about the smallest in absolute value integer $n_1(k)$ for which there exists  a triple $\{a, b, c\}$ of nonzero integers and a set $N$ of integers of size $k$ containing $n_1(k)$ such that $\{a, b, c\}$ is a $D(n)$-set for all $n\in N$?

\section{Diophantine $m$-tuples and elliptic curves}\label{sec2}

In this section we present several infinite families of Diophantine triples (i.e.\ $D(1)$-sets) which are also  $D(n)$-sets for two distinct $n$'s with $n\neq 1$. To that end, we first translate the problem to a problem about induced elliptic curves.

Let $\{a, b, c\}$ be a Diophantine triple and let $bc+1=r^2$, $ca+1=s^2$, $ab+1=t^2$.
Of interest to us are integer solutions $x$ of the system of equations
\begin{equation}\label{system}
\quad x+bc=\square, \quad x+ca=\square, \quad x+ab=\square.
\end{equation}
Note that one solution is $x=1$ (where the corresponding squares are $r^2, s^2, t^2$).
The following proposition (see e.g.\ \cite[4.1, p.~37]{hu}, \cite[4.2, p.~85]{kn})
reduces the study of solutions of this system of equations to the study of the induced elliptic curve \eqref{curve}.

\begin{proposition}\label{ell}
For $T\in E(\Q)$ we have that $x=x(T)$ is a rational solution of \eqref{system}
if and only if $T\in 2E(\Q)$.
\end{proposition}


Note that $E(\Q)$ has several obvious rational points:
\[
A=(-bc, 0), \quad B=(-ca, 0), \quad C=(-ab, 0),\quad
P=(0, abc), \quad S=(1, rst).
\]

By Proposition~\ref{ell} it follows that for every point $T\in 2E(\Q)\cap \Z^2$, we have that $\{a, b, c\}$ is a  $D(x(T))$-set. Since we are assuming that $\{a, b, c\}$ is a $D(1)$-set it follows that $S\in 2E(\Q)\cap \Z^2$. Indeed, it can be shown that $S=2R$ where
\[
R=(rs+rt+st+1, (r+s)(r+t)(s+t))\in E(\Q)\cap \Z^2.
\]

Note that $A, B, C$ are points of order $2$. In general, we may expect that the points $P$ and $S$ are two
independent points of infinite orders, and therefore that the rank of the induced elliptic curve is $\geq 2$. However, if $c = a+b \pm 2r$, where $ab+1=r^2$, then the direct computation shows that $2P=\pm S$. We remark that a Diophantine triple $\{a, b, c\}$ such that $c = a+b \pm 2r$ is said to be \emph{regular}.



In view of the above observations, we are led to look for triples $\{a, b,c\}$ for which  $2kP+\ell S\in \Z^2$ for some $k,\ell \in \Z$. One easily finds that
\[
2P=\left (\frac{1}{4}(a + b + c)^2 - ab - ac - bc , -abc - \frac{1}{8}((a + b + c)^2 - 4ab - 4ac - 4bc)(a + b + c)\right).
\]
Thus, the following holds.
\begin{lemma}\label{l:2P}
Let $a, b, c$ be nonzero integers such that $a+b+c$ is even. Then $\{a, b, c\}$ is a $D(n)$-set for
\begin{equation}\label{2P}
n=\frac{1}{4} (a + b + c)^2 - ab - ac - bc,
\end{equation}
provided $n\neq 0$. Furthermore, $n=0$ is equivalent to $c=a+b\pm 2\sqrt{ab}=(\sqrt{a}\pm\sqrt{b})^2$.

\end{lemma}

A direct proof of Lemma~\ref{l:2P} can be found in \cite{zg15}. Note that $n=1$ in Lemma~\ref{l:2P} if and only if $c=a+b\pm 2\sqrt{ab+1}$. Furthermore, note that $n\neq 0$ for any Diophantine triple $\{a, b, c\}$. Thus, we have the following corollary.

\begin{corollary}\label{c:abceven}
Any Diophantine triple $\{a, b, c\}$ such that $a+b+c$ is even and $c\neq a+b\pm 2\sqrt{ab+1}$
is also a $D(n)$-set for some $n\neq 1$.
\end{corollary}


A computer search shows that when $a$ and $b$ are in range $1$ to $1000$, $c$ in range $1$ to $1000 000$,
and $\{a, b, c\}$ is a $D(1)$-set, the corresponding points $S-2P$ and $4P$ never have integer coordinates. On the other hand,
the point $S+2P=2(R+P)$
has integer coordinates for the following $(a,b,c)$ such that $a<b<c$, $\{a, b, c\}$ is a $D(1)$-set, $a$ and $b$ in range from $1$ to $1000$, and $c$ in range $1$ to $1000 000$:
\begin{align*}
&(4, 12, 420),
(4, 420, 14280),
(12, 24, 2380),
(12, 420, 41184),
(24, 40, 7812),\\
&(40, 60, 19404),
(60, 84, 40612),
(84, 112, 75660),
(112, 144, 129540),\\
&(144, 180, 208012),
(180, 220, 317604),
(220, 264, 465612),
(264, 312, 660100),\\
&(312, 364, 909900).
\end{align*}
These examples suggest that there are infinitely many Diophantine triple $\{a,b,c\}$ such that the corresponding point
$S+2P$ has integer coordinates. We now show that this is indeed true. The point $S+2P$ has the following $x$-coordinate:
\begin{align*}
&-\frac{1}{4}(a + b + c)^2 -1+\\
&+ \frac{1}{4}\left(8abc + ((a + b + c)^2 - 4ab - 4ac - 4bc)(a + b + c) + 8\sqrt{ab + 1}\sqrt{ac + 1}\sqrt{bc + 1}\right)^2/ \\
&/ ((a + b + c)^2 - 4ab - 4ac - 4bc - 4)^2.
\end{align*}
We first note that all the examples above satisfy an additional condition that
$x(S+2P)=a+b+c$. A straightforward calculation shows that the condition $x(S + 2P)=a+b+c$ is equivalent to $q_1q_2q_3=0$,
where
\begin{align*}
q_1&=-4+a^2-2ab+b^2-2ac-2bc+c^2, \\
q_2&=a^2-4a-2ac-4c+c^2-2ab-4b-8abc-2bc+b^2, \\
q_3&=-4a-4b-4c-2ab-2ac-2bc-4abc+a^2+b^2+c^2-2a^2b-2a^2c-2ab^2-2ac^2\\
 &\mbox{}-2b^2c-2bc^2-2a^2b^2+2a^3+2b^3
+2c^3+a^4+b^4+c^4-2a^2c^2-2b^2c^2.
\end{align*}

The condition $q_1=0$ is equivalent to $c=a+b\pm 2\sqrt{ab+1}$. Thus, if $c=a+b\pm 2\sqrt{ab+1}$, then $\{a, b, c\}$ is a $D(n)$-set for $n=a+b+c=x(S+2P)$. However, by  Corollary~\ref{c:abceven}, such $\{a, b, c\}$ is not a $D(n)$-set for $n=x(2P)$, so in this way we can not get a Diophantine triple which is  also a $D(n)$-set for two distinct $n$'s with $n\neq 1$.
The solutions $(a, b, c)$ of $q_3=0$ are not Diophantine triples.
To see that, let $\sigma_1 = a+b+c, \sigma_2 = ab+bc+ca, \sigma_3 = abc$. Then $q_3 = 0$ is equivalent to $\sigma_1^4+2\sigma_1^3+(1-4\sigma_2)\sigma_1^2+(8\sigma_3-8\sigma_2-4)\sigma_1 + 8\sigma_3 - 4\sigma_2 = 0$, which can be written as $(\sigma_1^2-4\sigma_2) (\sigma_1 +1)^2+(8\sigma_3-4)(\sigma_1+1)+4= 0$. Therefore, $(\sigma_1 +1) | 4$, i.e.~$\sigma_1 +1 \in \{-4, -2, -1, 1, 2, 4\} \Rightarrow \sigma_1=a+b+c \in \{-5, -3, -2, 0, 1, 3\}$. If we assume that $(a, b, c)$ is a Diophantine triple, then $a, b$ and $c$ must have the same sign. Then one can easily check that $a+b+c \in \{-5, -3, -2, 0, 1, 3\}$ does not result in a Diophantine triple (in fact,  for Diophantine triples $\{a, b, c\}$ we have $|a+b+c| \geqslant 12$).

On the other hand, the condition $q_2=0$ is equivalent to
\begin{equation}\label{regc}
c=2+a+b+4ab\pm 2\sqrt{(2a+1)(2b+1)(ab+1)},
\end{equation}
and this is exactly the condition that
$\{2,a,b,c\}$ is a regular Diophantine quadruple. Like we mentioned in the introduction, it is a conjecture that all Diophantine quadruples are regular.
\begin{theorem}\label{th:reg}
Let $\{2,a,b,c\}$ be a regular Diophantine quadruple.
Then the Diophantine triple $\{a,b,c\}$ is also a $D(n)$-set for two distinct $n$'s with $n\neq 1$.
\end{theorem}

\proof
We claim that $n_2=x(S+2P)$ and $n_3=x(2P)$ (corresponding to the Diophantine triple $\{a, b, c\}$) have the required properties. Note that we must have $a, b, c>2$. In what follows, we assume that $a<b<c$. We first show that $\{a,b,c\}$ is also a $D(n_3)$-set. Recall that
\begin{equation}\label{2P}
n_3=\frac{1}{4} (a + b + c)^2 - ab - ac - bc,
\end{equation}
Since $2a+1$, $2b+1$ and $2c+1$ are perfect squares, we conclude that $a$, $b$ and $c$ are even.
By \cite[Lemma 14]{djnt}, we have
$c>4ab$ and therefore
$c\neq a+b\pm 2\sqrt{ab+1}$.  By Corollary~\ref{c:abceven} it follows that
$\{a,b,c\}$ is a $D(n_3)$-set with $n_3\neq 1$.

We now show that $\{a,b,c\}$ is also a $D(n_2)$-set.
Since $\{2,a,b,c\}$ is a regular Diophantine quadruple,  it follows that $n_2=a+b+c$. (See the discussion just before this theorem). By $a, b, c>2$ it follows that  $n_2\neq 1$. Thus, $\{a,b,c\}$ is a $D(n_2)$-set with $n_2\neq 1$.

It remains to show that $n_2\neq n_3$. Note that $n_2=n_3$ is equivalent to $2P=\pm (S+2P)$, i.e. $2(R+2P)=\mathcal{O}$.
This condition leads to
\begin{equation}\label{eq:rst14}
(r\pm s)(r\pm t)=\frac{1}{4}(c-a-b)^2.
\end{equation}
Rough estimates show that the left hand side is less than $4c\sqrt{ab}<2c\sqrt{c}$,
while the right hand side is greater than $c^2/16$, so that for $c\geq 1024$ we immediately get a contradiction.
The only Diophantine triple with $c<1024$ such that $\{2,a,b,c\}$ is a regular Diophantine quadruple is $\{4,12,420\}$, for which we directly check that
\eqref{eq:rst14} does not hold.
\qed

We now give several explicit infinite families of Diophantine triples $\{a,b, c\}$ satisfying the conditions of Theorem \ref{th:reg}. We first choose $a$ and $b$ so that $\{2,a, b\}$ is a Diophantine triple. We obtain one family of Diophantine triples $\{2,a, b\}$ by letting
 $b=2+a+2\sqrt{2a+1}$ and $2a+1=(2i+1)^2$ for $i\in \N$. Then $a=2(i + 1)$ and $b=2(i + 2)(i + 1)$. We compute $c$ using the regularity condition \eqref{regc} and obtain the following corollary of Theorem~\ref{th:reg}.

\begin{corollary}\label{inf1}
Let $i$ be a positive integer and let
\[
a=2(i + 1)i, \quad b=2(i + 2)(i + 1), \quad c=4(2i^2 + 4i + 1)(2i + 3)(2i + 1).
\]
Then $\{a, b, c\}$ is a $D(n)$-set for $n=n_1, n_2, n_3$, where
\begin{align*}
n_1&=1,\\
n_2&=32i^4 + 128i^3 + 172i^2 + 88i + 16,\\
n_3&=256i^8 + 2048i^7 + 6720i^6 + 11648i^5 + 11456i^4 + 6400i^3 + 1932i^2 + 280i + 16.
\end{align*}
\end{corollary}

Another family of Diophantine triples of type $\{2, a, b\}$ can be obtained by taking $a=2(i + 1)i$, $b=4(2i^2 + 4i + 1)(2i + 3)(2i + 1)$.  We again compute $c$ using the regularity condition \eqref{regc} and obtain the following corollary of Theorem~\ref{th:reg}.

\begin{corollary}\label{inf2}
Let $i$ be a positive integer and let
\[
a=2(i + 1)i, \quad b= 4(2i^2 + 4i + 1)(2i + 3)(2i + 1), \quad c=2(4i+1)(4i+3)(4i^2+9i+4)(4i^2+7i+1).
\]
Then $\{a, b, c\}$ is a $D(n)$-set for $n=n_1, n_2, n_3$, where
\begin{align*}
n_1&=1,\\
n_2&=512i^6+2560i^5+4832i^4+4352i^3+1980i^2+432i+36,\\
n_3&=65536i^{12}+655360i^{11}+2859008i^{10}+7151616i^9+11346176i^8+11932672i^7+8450112i^6 \\
 &\mbox{}+4012672i^5+1249280i^4+243840i^3+27612i^2+1584i+36.
\end{align*}
\end{corollary}

We now present another way to get an infinite family of triples $\{a,b,c\}$ satisfying the conditions of Theorem \ref{th:reg}. We first choose $a$ so that $\{2, a\}$ is a Diophantine pair. If $2a+1=k^2$, consider the sequence $(b_i)_{i\geq 0}$ defined by the recurrence relation
\begin{equation}\label{rec}
 b_0 = 0, \quad b_1=2+a+2k, \quad b_2=4k(k+2)(k+a), \quad b_{i+3}=(4k^2-1)b_{i+2} - (4k^2-1)b_{i+1} + b_i.
 \end{equation}
Then $\{2,a,b_i,b_{i+1}\}$ is a regular Diophantine quadruple for all $i\geq 1$, see e.g. \cite[Section 3]{dacta}, \cite[Lemma 15]{djnt}.
For example, take $a=4$. Then $k=3$ and
\[
 b_0 = 0, \quad b_1 = 12, \quad b_2 = 420, \quad b_{i+3} = 35b_{i+2}-35b_{i+1}+b_i, \quad i\geq 3,
 \]
so $\{2,4,b_i,b_{i+1}\}$ is a regular Diophantine quadruple for all $i\geq 1$. In fact, by \cite{fuj} all Diophantine quadruples containing $2$ and $4$ are of this form.

\begin{corollary}\label{inf3}

Let the sequence $(b_i)_{i\geq 0}$ be defined by \eqref{rec}. Then for all positive integers $i$
the triple $\{4, b_i, b_{i+1}\}$ is a $D(n)$-set for $n=n_1, n_2, n_3$, where
\begin{align*}
n_1&=1,\\
n_2&=4+b_i+b_{i+1},\\
n_3&=\frac{1}{4}(4+b_i+b_{i+1})^2 - 4b_i - 4b_{i+1}-b_ib_{i+1} .
\end{align*}
\end{corollary}

\section{$D(n)$-sets for $n\in \{1,n_2,n_3,n_4\}$}\label{sec3}

We now discuss the existence of Diophantine triples $\{a,b,c\}$ which are also
$D(n)$-sets for three distinct $n$'s with $n\neq 1$. We have been able to find several such
examples, but we do not know whether there exists infinitely many such triples.  We performed an extensive computer search to find examples of Diophantine triples from Corollaries \ref{inf1}, \ref{inf2} and \ref{inf3} for which there exists an additional integer $n_4$
such that these triples are also $D(n_4)$-sets.

A computer search shows that for $i<100000$, there is no $n\neq n_1, n_2, n_3$ such that
\[
\{a, b, c\}=\{2(i + 1)i,  2(i + 2)(i + 1), 4(2i^2 + 4i + 1)(2i + 3)(2i + 1)\}
\]
 is a $D(n)$-set,
except for $i=1$ and $i=4$. When $i=1$, then also $n_4=3796$ is such that $\{a, b, c\}$ is a $D(n_4)$-set, and when $i=4$ then also $n_4=3680161$ is such that $\{a, b, c\}$ is a $D(n_4)$-set.
Thus, the triple $\{4, 12, 420\}$ is a $D(n)$-set for $n=1, 436, 3796, 40756$
and the triple $\{40, 60, 19404\}$ is  a $D(n)$-set for $n=1, 19504, 3680161, 93158704$.
This suggests that it is unlikely that we can find an infinite subfamily of the family of triples from Corollary~\ref{inf1}, such that all triples from that family are $D(n)$-sets for four distinct $n$'s.

\begin{remark}
We now present one way to test that  for $i<100 000$, there is no $n\neq n_1, n_2, n_3$
such that $\{a, b, c\}=\{2(i + 1)i,  2(i + 2)(i + 1), 4(2i^2 + 4i + 1)(2i + 3)(2i + 1)\}$
is a $D(n)$-set except for $i=1$ and $i=4$.
Assume that $\{a, b, c\}$ is a $D(n)$-triple for some $n \notin \{n_1, n_2, n_3\}$. Let
\[
bc+n=r^2,\quad ca+n=s^2, \quad ab+n=t^2.
\]
Then
\[
r-t\mid b(c-a)=4(4i^2 + 9i + 3)(4i^2 + 7i + 2)(i + 2)(i+ 1).
\]
 From $d=r-t$ we can find $n$. Indeed $2r=(r-t)+(r+t)=d+b(c-a)/d$ and $n=r^2-bc$. We thus let $i$ run from $1$ to $100 000$ and let $d$ run over all divisors of
\[
N=b(c-a)=4(4i^2 + 9i + 3)(4i^2 + 7i + 2)(i + 2)(i+ 1)
\]
and then test if $n=(N/d+d)^2/4-bc$ is such that $ca+n\in \Z^2$. Note that $ab+n$ will always be a square since $ab+n=(N/d-d)^2/4$.
\end{remark}

It is natural to look for triples which are $D(n)$-sets for additional $n$'s within the family from Corollary~\ref{inf1} by considering
those whose induced elliptic curve has larger rank.
We first consider the elliptic curve over $\mathbb{Q}(t)$ defined by
\[
E(\Q(t))\colon\quad y^2=(x+ab)(x+ac)(x+bc),
\]
where
\begin{align*}
a=2(t + 1)t,\quad b=2(t + 2)(t + 1), \quad c=(2t^2 + 4t + 1)(2t + 3)(2t + 1).
\end{align*}
We now apply \cite[Theorem 1.1]{GT} to show that the rank of this curve (over $\Q(t)$) is  $2$.
For the specialization $t\mapsto 22$ the conditions of the main theorem of \cite{GT} are satisfied,
so the corresponding specialization homomorphism is injective.
The specialized curve $E_{22}$ over $\mathbb Q$ is given by
\[
y^2 = x^3+1182678423x^2+349062404795541360x+24368563262176237011600
\]
and has rank $2$, where the $x$-coordinates of the two free generators $G_1$ and $G_2$ are
$x(G_1)=0$, $x(G_2)=603593508$. We then note that $G_1$ and $G_2$ are the specialized points for the points $P$ and $R$ on $E(\mathbb Q(t))$ (where, as always, $P=(0, abc)$, $S=(1, \sqrt{(bc+1)(ca+1)(ab+1)})$ and $R$ is defined by $S=2R$). Thus, $P$ and $R$ are the free generators of  $E(\mathbb Q(t))$.

Now we show that for Diophantine triples  from Corollary~\ref{inf1} such that $i(3i+1)$ is a square (i.e.~for $i=1,16,225$, etc.) we have that the rank of the induced elliptic curve $E(\Q)$ is $\geq 3$. To show that, we transform $E(\Q)$ into $y^2 = x(x+ac-ab)(x+bc-ab)$ (by substituting $x$ by $x-ab$, so that the point $(0,0)$ is of order $2$) and search for additional points by considering the factors of
\[
ab(c-a)(c-b)=16i(i+1)^2(i+2)(4i^2+7i+1)(4i^2+9i+4)(4i^2+9i+3)(4i^2+7i+2).
\]
Now, for the factor $x=-4(4i^2+9i+4)(4i^2+9i+3)(i+1)^2$ we have that
\[
 x(x+ac-ab)(x+bc-ab)=64(4i^2+9i+4)^2(4i^2+9i+3)^2(i+1)^4(2i+3)^2(3i+1)i,
 \]
and this is a square if and only if $i(3i+1)$ is a square. So, if $i(3i+1)$ is a square we have an additional point on $E(\Q)$ (whose $x$-coordinate is $-4(4i^2+9i+4)(4i^2+9i+3)(i+1)^2$). If we denote that additional point by $Q$, then the point $2Q$ might be a candidate for an additional $n$.
However, $x(2Q)=(64i^8+576i^7+2100i^6+4020i^5+4389i^4+2794i^3+1025i^2+200i+16)/(i(3i+1))$ is an integer
only for $i=1$, which gives $n_4=3796$ for the triple $\{4, 12,  420\}$. It seems that this and
similar approaches with larger rank curves can produce only some sporadic examples.

Recall that the example $\{4, 12,  420\}$ also belongs to the family of triples from Corollary
\ref{inf3}. We have been able to find two more examples of triples which are $D(n)$-sets for one additional $n$ within that family:
$\{4, 420,  14280\}$ is a $D(n)$-set for $n=1,14704, 950896, 47995504$, while
$\{4, 485112, 16479540\}$ is a $D(n)$-set for $n=1,16964656,2007609136,63955397832496$.

We have also searched for Diophantine triples $\{a,b,c\}$ which are $D(n)$-sets for four distinct $n$'s without requiring that $\{2,a,b,c\}$ is a regular Diophantine quadruple. It is natural to look at families of non-regular Diophantine triples $\{a,b,c\}$, such that the points $P$ and $S$ are independent, which additionally satisfy that $a+b+c$ is even, so that the Corollary \ref{c:abceven} can be applied. In particular, we have performed an extensive search of the family $\{k \pm 1, 4k, 16k^3-4k\}$ with odd $k$, which is in some way the simplest such family.
Our search has resulted in two additional examples of Diophantine triples which are $D(n)$-sets for four distinct $n$'s. We have that $\{10, 44,  21252\}$ is a $D(n)$-set for $n=1, 825841, 6921721, 112338361$ and $\{78, 308, 7304220\}$ is a $D(n)$-set for $242805865, 4770226465,13336497750865$.

In conclusion, we have been able to find the following triples $\{a,b,c\}$  which are $D(n)$-sets for $n_1=1<n_2<n_3<n_4$:

\begin{center}
\begin{tabular}{l|l}
$\{a,b,c\}$ & $n_2,n_3,n_4$ \\
\hline
$\{4, 12,  420\}$ &      $436, 3796, 40756$ \\
$\{10, 44,  21252\}$ &  $825841, 6921721, 112338361$ \\
$\{4, 420,  14280\}$ &   $14704, 950896, 47995504$ \\
$\{40, 60, 19404\}$ &    $19504, 3680161, 93158704$ \\
$\{78, 308, 7304220\}$ & $242805865, 4770226465,13336497750865$ \\
$\{4, 485112, 16479540\}$ &   $16964656, 2007609136, 63955397832496$
\end{tabular}
\end{center}

\section{A modification of the problem}\label{sec4}

So far we were interested in the maximum size of a set $N$ of nonzero integers containing $1$ for which there exists a triple of nonzero integers
$\{a, b, c\}$ which is a $D(n)$-set for all $n\in N$.
Now note that if we omit the condition $1\in N$,
then the size of a set $N$ for which there exists a triple $\{a, b, c\}$ of nonzero integers
which is a $D(n)$-set for all $n\in N$ can be arbitrarily large. Indeed,
take any triple $\{a, b, c\}$ such that the induced elliptic curve $E(\Q)$
has positive rank. Then there are infinitely many rational points
on $E(\Q)$. For an arbitrary large positive integer $m$ we may choose $m$ distinct rational points $R_1,\ldots,R_m\in 2E(\Q)$, so that by Proposition~\ref{ell} we have
\[
x(R_i)+bc=\square,\quad  x(R_i)+ca=\square, \quad x(R_i)+ab =\square,
\]
where $\square$ stands for a square of a rational number. We do so by taking points of the form $2m_1P_1+2m_2P_2+\cdots+2m_rP_r$, where $P_1, \ldots, P_r$ are the generators of $E(\Q)$.
We then let $z\in \Z\setminus\{0\}$ be such that $z^2x(R_i)\in \Z$ for all $i=1, 2,\ldots, m$.
Then the triple $\{az, bz, cz\}$ is a $D(n)$-set for
$n=x(R_i)z^2$ for all $i=1,2,\ldots,m$.

We consider the following modification of the original problem.
\begin{question}\label{Q1}
For a given positive integer $k$, what can be said about the smallest in absolute value nonzero integer $n_1(k)$ for which there exists  a triple $\{a, b, c\}$ of nonzero integers and a set $N$ of integers of size $k$ containing $n_1(k)$ such that $\{a, b, c\}$ is a $D(n)$-set for all $n\in N$?
\end{question}

Note that if $k\leq 4$, then $n_1(k)=1$ since we have found examples of  Diophantine triples $\{a, b, c\}$ which are also $D(n)$-sets for three distinct $n$'s greater than $1$. We have not been able to find neither $D(1)$-set nor $D(-1)$-set $\{a, b, c\}$ which is also a $D(n)$-set for four distinct $n$'s with $|n|>1$, and in fact we suspect that  $|n_1(5)|>1$ based on exhaustive but unsuccessful computer search.
(We remark that it has been conjectured that there do not exist $D(-1)$-quadruples, and that it is known that there do not exist $D(-1)$-quintuples and that there are at most finitely many $D(-1)$-quadruples~\cite{ dff07,df05}).

Our initial motivation for Question~\ref{Q1} was to see how close to $1$ we can get with $n_1(5)$.
We can show that $|n_1(5)|\leq 36$. To that end we consider again the family of Diophantine triples
$\{k-1, 4k, 16k^3-4k\}$ with $k\geq 2$.
 For $k=2$, we get a triple $\{1,8,120\}$ whose induced elliptic curve $E(\Q)$
has rank $3$. Following the procedure described above we find points $R_1, \ldots, R_5\in 2E(\Q)$ such that
\[
x(R_1)=1, \quad x(R_2)=721, \quad x(R_3)=12289/4, \quad x(R_4)=769/9, \quad x(R_5)=1921/36.
\]
We then let $z=6$. It follows that the triple $\{az, bz, cz\}=\{6,48,720\}$ is a $D(n)$-set for
$n=36, 1921, 3076, 25956, 110601$. (We chose $R_2, \ldots, R_5\in 2E(\Q)$ so that their $x$-coordinates have relatively small denominators. We obtained the $n$'s using $n=x(R_i)z^2$, $i=1,2,\ldots,5$).

To bound $|n_1(k)|$ for larger values of $k$ we use a similar idea. We do not require that the starting triple $\{a, b, c\}$ of nonzero integers
 is a $D(1)$-triple, but  search for triples of nonzero integers $\{a, b, c\}$ whose induced elliptic curve has
relatively large rank $r$ (say $r \geq 5$). Namely, roughly speaking, with more generators $P_1, \ldots, P_r$, we may expect
more linear combinations $2m_1P_1+2m_2P_2+\cdots+2m_rP_r$ with reasonably small height.
The usual arguments
(see e.g. \cite{aabl} or \cite{dels}) suggest that $|n_1(k)|$ grows as
$k^{(r+2)/r}$. However, this estimate ignores the heights of generators, which conjecturally
grow with rank ($\max_{1\leq i \leq r} h(P_i) \gg r\log(r)$, see e.g.~\cite{sil}).
In conclusion, we don't expect that the optimal values of $n_1(k)$ will be obtained through a triple whose induced elliptic curve has maximal known rank
(which is $r=15$ for $\{a,b,c\} = \{153437625700011715932305091$, $148276675335745173559912875,
 153776032705386698956099516\}$, see \cite{elk}), but through triples whose induced elliptic curve has moderate rank and generators of relatively small heights.

 The following table summarizes our findings:
\begin{center}
\begin{tabular}{ |c|c|c|c| }
 \hline
$k$ & $|n_1(k)| \leq $ & rank & $\{a,b,c\}$ \\
\hline
$5$ & $36$ & $3$ & $\{6, 48, 720\}$  \\
$6$ & $215$ & $3$ & $\{28, 168, 1848\}$ \\
$7$ & $900$ & $4$ & $\{380, 1400, 3240\}$ \\
$8$ & $7740$ & $3$ & $\{168, 1008, 11088\}$ \\
$9$ & $32400$ & $4$ & $\{2280, 8400, 19440\}$ \\
$10$ & $129600$ & $4$ & $\{4560, 16800, 38880\}$ \\
$11$ & $215991$ & $5$ & $\{9120, 22770, 30960\}$ \\
$12$ & $863964$ & $5$ & $\{18240, 45540, 61920\}$ \\
$13$ & $4932144$ & $5$ & $\{37128, 118440, 182280\}$ \\
$14$ & $7706475$ &  $5$ & $\{46410, 148050, 227850\}$ \\
$15$ & $30825900$ & $5$ & $\{92820, 296100, 455700\}$ \\
$16$ & $123303600$ & $5$ & $\{185640, 592200, 911400\}$ \\
$17$ & $371289600$ & $5$ & $\{59400, 108360, 223200\}$ \\
$18$ & $4438929600$ & $5$ & $\{1113840, 3553200, 5468400\}$ \\
$19$ & $18193190400$ & $5$ & $\{415800, 758520, 1562400\}$ \\
$20$ & $18193190400$ & $5$ & $\{415800, 758520, 1562400\}$ \\
\hline
\end{tabular}
\end{center}


 \vspace{0.5cm}
{\bf Acknowledgements.}
N.A.\@, A.D. and P.T.\@ were supported by the Croatian Science Foundation under the project no.~6422. A.D.\@ acknowledges support from the QuantiXLie Center of Excellence. D.K.\@ acknowledges support of the Austrian Science Fund (FWF):
J-3955, and support of the Austrian Academy of Sciences through
JESH project.

\end{document}